\documentclass{gtart}
\input gtmonout
\volumenumber{1}
\volumeyear{1998}
\volumename{The Epstein birthday schrift}
\pagenumbers{249}{260}
\published{26 October 1998}
\papernumber{11}
\received{26 August 1997}

\newcommand{\Z}{\mbox{${\rm Z\!\! Z}$}}    % die ganzen Zahlen
\def \fsa {{\sf fsa}}
\def \kbmag {{\sf kbmag}}

\def\ssp{\stdspace}

\newtheorem{theorem}{Theorem}[section]

\theoremstyle{definition}
\newtheorem{defn}[theorem]{Definition}

\begin{document}

\title{Automatic groups, subgroups and cosets}

\author{Derek F Holt}

\address{Mathematics Institute, University of Warwick\\
Coventry, CV4 7AL, UK.}

\email{dfh@maths.warwick.ac.uk}

\begin{abstract} 
The history, definition and principal properties of automatic groups
and their generalisations to subgroups and cosets are reviewed
briefly, mainly from a computational perspective. A result about the
asynchronous automaticity of an HNN extension is then proved and
applied to an example that was proposed by Mark Sapir.
\end{abstract}



\maketitle

The concept of an automatic group was introduced in 1986 by Thurston,
motivated by some results of Jim Cannon on hyperbolic groups.
Much of the basic theory of this important class of groups was
developed by David Epstein during the following few years.

In the first section of this paper, we review briefly the history,
definition and properties of automatic groups and their generalisation
to subgroups and cosets, mainly from a perspective of carrying out
efficient computations within such groups and their subgroups.
In the second section, we prove a result about the (asynchronous) automaticity
of an HNN extension, and use it, together with the results of some machine
computations, to prove that a particular group, defined by Mark
Sapir, is asynchronously automatic.

\section{Definitions and discussion}
\subsection{Automatic groups}
In~\cite{Can}, J.W. Cannon proved certain geometrical properties of
the Cayley graph of cocompact discrete hyperbolic groups.
Two years later, in 1986, W. Thurston noticed that some of these
properties could be reformulated in terms of finite state automata (\fsa; this
abbreviation will be used for both the singular and plural).

In particular, the goedesic paths in the Cayley graph that start at the
origin form a regular set or, equivalently, they form the language of
an \fsa. Furthermore, any pair of such geodesic paths that
end at the same or neighbouring vertices lie within a bounded distance of
each other. It can be deduced that such geodesic pairs also form
the language of an \fsa.
This led Thurston to formulate the following general definition.

\begin{defn} \label{aut}
Let $G$ be a group with finite generating set $X$, let $A = X \cup X^{-1}$,
and let $A^\prime = A \cup \{\$\}$, where $\$ \not\in A$.
Then $G$ is said to be {\em automatic} (with respect to $X$),
if there exist {\sf fsa} $W$ and $M_a$ for each $a \in A^\prime$, such that\\
(i)\ssp $W$ has input alphabet $A$, and accepts at least one word in $A^*$ mapping
onto each element of $G$.\\
(ii)\ssp Each $M_a$ has input alphabet $A^\prime \times A^\prime$,
it accepts only padded pairs, and it accepts
the padded pair $(w^+,x^+)$ for $w,x \in A^*$ if and only if 
$w,x \in L(W)$ and $wa =_G x$.
\end{defn}

Here $A^*$ as usual denotes the set of words in $A$.
For $w \in A^*$, ${\overline {w}}$ denotes the element of $G$ onto which
$w$ maps; for $w, x \in A^*$, we also use $w =_G x$ to mean that
$w,x$ map onto the same element of $G$.
The extra symbol $\$$ maps onto the identity element of $G$.
For $w,x \in A^*$, the associated {\em padded pair} $(w^+,x^+) \in
(A^\prime \times A^\prime)^*$ is obtained by adjoining symbols $\$$ to the
end of the shorter of $w$ and $x$ to make them have equal length.
The language of the \fsa\ $W$ is denoted by $L(W)$.
For general properties of finite state automata, the user is referred to
any textbook on automata or formal language theory, such as~\cite{HU}.

In the definition, $W$ is called the {\em word-acceptor} and the $M_a$
the {\em multiplier} automata.
The complete collection $\{W,M_a\}$ is known as an
{\em automatic structure} for $G$. 
Note that the multiplier $M_\$$ recognises equality in $G$ between words
in $L(W)$. From a given automatic structure, we can always use $M_\$$ to
construct another one such that $W$ accepts a unique word mapping onto each
element of $G$; we simply choose the lexicographically least amongst the
shortest words that map onto each element as the `normal form' representative
of that element. We shall call such a $W$ a word-acceptor with uniqueness.

The best general reference for the theory of automatic groups is the
multi-author book~\cite{ECHLPT}. In particular, it turns out that the
automaticity of $G$ is independent of the choice of generating set $X$.
This immediately suggests that the definition is a sensible one, because
it means that automaticity is an algebraic property of the
group, rather than just a geometrical property of its Cayley graph.

All finite groups are easily seen to be automatic; in fact the class of
automatic groups is invariant under finite variations, such as sub-
and super-groups of finite index. It is also closed under direct and
free products, and includes, for example, all word-hyperbolic groups,
braid groups, Coxeter groups and Artin groups of finite and of `large' type.
All automatic groups have finite presentations.

Some of the most important and useful applications of this theory only involve
an explicit knowledge of a word acceptor with uniqueness, particularly
in the frequently occurring case when the accepted words are all
geodesics in the Cayley graph.
 From such a word-acceptor, one can quickly enumerate unique
representatives of all words up to a given length. This can 
serve as an invaluable time-saving device in certain computer
graphics applications, such as drawing tessellations of hyperbolic
space on which these groups act freely.
One can also use $W$ to compute the growth function for the group
(see~\cite{EIZ}).

Another important application of automatic structures for groups $G$ is their
use for the efficient (quadratic time) solution of the word problem in $G$.
More precisely, the multiplier automata can be used to reduce an
arbitrary word in $A^*$ in quadratic time to the $G$--equivalent word in $L(W)$.

With these applications in mind, a collection of programs was written at
Warwick in the late 1980's for computing automatic structures. These programs
take a finite presentation of the group $G$ as input. Currently, they only
work for so-called {\em shortlex} structures, which are those in which
$L(W)$ consists of the lexicographically least amongst the
shortest words that map onto each group element. (So $W$ depends upon
the order of $A$ as well as on $A$ itself.) Many, but not all, of the known
classes of automatic groups are known to possess shortlex structures.
The programs are described in some detail in~\cite{EHR} and~\cite{Holt},
and in a much more general setting in~\cite{ECHLPT}. The latest version is
part of a package called \kbmag\ and is available by anonymous
\verb-ftp- from \verb-ftp.maths.warwick.ac.uk- in the directory
\verb-people/dfh/kbmag2-.

 From an algorithmic point of view, there is a close connection between
automatic groups and rewriting systems for groups, and the programs
used make use of the Knuth--Bendix completion process in groups.
However, typically, this process alone would not terminate and in fact
automatic groups normally have infinite regular rather than
finite complete rewriting systems. When the automatic structure is
successfully computed it is, in some sense, enabling this infinite
regular system to be used to solve the word problem in a manner that is
typically at least as efficient as could be done with a finite rewriting
system. The idea of trying to use infinite regular rewriting systems
for this purpose was first proposed by Gilman in~\cite{Gil}.

Given a word-acceptor automaton for a group, it turns out that
the existence and properties of the multiplier automata
are equivalent to the so-called (synchronous) fellow-traveller property,
which was one of the geometrical properties of hyperbolic groups
observed originally by J. W. Cannon, and is defined as follows.

For a word $w \in A^*$ we denote the length of $w$ by $l(w)$ and, for
$g \in G$, $l(g)$ (or more precisely $l_A(g)$) denotes the length of the
shortest word $w \in A^*$ with ${\overline w} = g$. For $t \geq 0$, $w(t)$
denotes the prefix of $w$ of length $t$ when $t \leq l(w)$, and $w(t)=w$
for $t \geq l(w)$. The fellow-traveller property asserts that there exists
a constant $k$ such that, for all $w, x \in L(W)$ and $a \in A$ such
that $wa =_G x$, and all $t \geq 0$, we have $l_A( 
{\overline {w(t)}}^{-1} {\overline {x(t)}} ) \leq k$.
In other words, two travellers proceeding at the same speed along the
words $w$ and $x$ from the base point in the
Cayley graph of $G$  would always remain a bounded
distance away from each other.

The fellow-traveller property enables the multiplier automata $M_a$ to be
defined in a uniform manner (see Definition 2.3.3
of~\cite{ECHLPT}). Their state set is the
set of triples $(s_1,s_2,g)$, where $s_1,s_2$ are states of $W$, and
$g \in G$ with $l(g) \leq k$. The start state is $(s_0,s_0,1)$, where
$s_0$ is the start state of $W$. For $(a_1,a_2) \in  A \times A$,
there is a transition from $(s_1,s_2,g)$ to $(t_1,t_2,h)$
with label $(a_1,a_2)$ if and only if there are transitions $s_1
\rightarrow t_1$ and $s_2 \rightarrow t_2$ in $W$ with labels
$a_1$ and $a_2$, respectively, and if $a_1^{-1}ga_2 =_G h$.
The state $(s_1,s_2,g)$ is a success state of $M_a$ if and only
if $s_1$ and $s_2$ are success states of $W$, and $g =_G a$.
Thus the $M_a$ differ only in their accept states. (We have omitted a
technicality from this definition. To deal with the padding symbol, we
have to add an extra state to $W$ which is reached when $W$ is in
an accept state and the padding symbol is read.)
It is clear that the $M_a$ behave precisely according to Condition (ii) of
Definition~\ref{aut}. This method is used to construct the $M_a$ in the
programs mentioned above.

Note also that it follows from the fellow-traveller property that if $g$
is any fixed element of $G$ and $w,x \in L(W)$ with $wg =_G x$, then
$w$ and $x$ fellow-travel with constant at most $kl_A(g)$.

Finally, we must mention the weaker concept of an asynchronously automatic
group, because it will arise in the next section.
The definition is the same as before, except that the multiplier automata are
allowed to read their two input strings at different rates. More precisely,
rather than reading one symbol from each of the two input words at each
transition, they read a symbol from one of the two words only, where the
choice of which word to read is a function of the state of $M_a$. Of course,
when the end of one of the words is reached, the other word must be selected.
See Chapter~7 of~\cite{ECHLPT} for the formal definition.
Again there is a corresponding fellow-traveller property, in which
the imaginary travellers are allowed to move at different speeds.
See~\cite{ECHLPT} or Section 7, Part II of~\cite{BGSS} for details.

The word problem is still solvable for asynchronously automatic groups, but it
is unknown whether this can be done in polynomial time.
There are examples known, such as the Baumslag--Solitar
groups $\langle x,y \, | \, y^{-1}x^py = x^q \, \rangle$ with $p \neq q$,
which are asynchronously automatic but not automatic.

There is a more detailed treatment, with references to the literature,
of the synchronous and asynchronous fellow-traveller properties in groups
in the article~\cite{Rees} in these proceedings.

\subsection{Subgroups}
Let $L=L(W)$ be the language of the word-acceptor in an automatic structure
of a group $G$. A subgroup $H$ of $G$ is called $L$--rational if $L \cap H$
is a regular language (ie the language of an \fsa). Such subgroups were
studied in~\cite{GS}, where it is proved that $L$--rational is equivalent to
$L$--quasiconvex. This means that any prefix of a word in $L \cap H$ lies
within a bounded distance of $H$ in the Cayley graph of $G$. Such subgroups
are always finitely generated.

An algorithm for constructing an \fsa\ $W_H$ with language $L \cap H$,
which takes as input an automatic structure for $G$ and
a set of generators for an $L$--rational subgroup $H$ of $G$, is
described in~\cite{Kap}. A practical and efficient version is
described in~\cite{Hurt}, and an implementation is available in \kbmag.

The \fsa\ $W_H$ can be used together with the automatic structure to
determine whether a given word in $A^*$ lies in $H$; that is, to solve
the generalised word problem for $H$ in $G$. First use the the multiplier
automata to reduce the word to one in $L$, and then use $W_H$ to test
whether it lies in $H$. Given $W_H$ and $W_K$ for two subgroups $H$ and $K$
of $G$, it is easy to intersect their languages to obtain a
\fsa\ $W_{H \cap K}$ for their intersection, which can then be used
to construct a finite generating set for $H \cap K$.

\subsection{Cosets}
It is possible to generalise the concept of an automatic group from a notion
about the elements of the group to one about the cosets of a given subgroup
$H$ of $G$. This has been carried out by two doctoral students of the author
(see~\cite{Redfern} and~\cite{Hurt}). The definition is as follows.

\begin{defn}
Let $G$ be a group with finite generating set $X$, let $A = X \cup X^{-1}$,
$A^\prime = A \cup \{\$\}$, and let $H$ be a subgroup of $G$.
Then $G$ is said to be
{\em coset automatic} with respect to $H$, if there exist {\sf fsa} $W$,
and $M_a$ for each $a \in A^\prime$, such that:\\
(i)\ssp $W$ has input alphabet $A$, and accepts at least one word in each right
coset of $H$ in  $G$;\\
(ii)\ssp Each $M_a$ has input alphabet $A^\prime \times A^\prime$,
it accepts only padded pairs, and it accepts
the padded pair $(w^+,x^+)$ for $w,x \in A^*$ if and only if 
$w,x \in L(W)$ and ${H\overline {wa}} = {H\bar x}$.
\end{defn}

Here $W$ is called the {\em coset word-acceptor} and the $M_a$
the {\em coset multiplier} automata.
The complete collection $\{W,M_a\}$ is known as an
{\em automatic coset system} for the pair $(G,H)$. Again the existence of such
a system turns out to be independent of the generating set $X$ of $G$,
and we can, if we wish, always find a new system in which $W$ accepts
a unique word in each right coset.

It is proved in~\cite{Redfern} that if $L$ is the language of the
shortlex automatic structure of a word-hyperbolic group $G$ (or even the
set of all geodesics in the Cayley graph of $G$), and if the subgroup
$H$ is $L$--quasiconvex, then $G$ is coset automatic with respect to $H$.
In~\cite{Hurt} the converse is proved for word-hyperbolic groups,
although we shall see from the example in the next section that the converse
does not hold in general.

An interesting application to the drawing of limit sets of Kleinian groups is
described in~\cite{MPR}. As in the graphical applications of ordinary
automatic structures, this involves only the use of $W$ to enumerate unique
shortest words in each coset.

An algorithm for computing automatic coset systems in the shortlex case
was first described in~\cite{Redfern}, and was implemented
by him as a standalone program.
It has the disadvantage that it is not usually possible to prove
conclusively that the system computed is correct.
A different approach is described in~\cite{Hurt}.
This does enable the output to be proved correct, but it requires an additional
hypothesis, to be described below, for it to work at all.
It has the further advantage that it has an optional extension to compute a
finite presentation for the subgroup $H$ of $G$ after the automatic coset
system has been found.
This second algorithm, together with the subgroup presentation facility,
has been implemented and is available in \kbmag.
The theory, implementation details and performance statistics
can also be found in~\cite{HH}.

These algorithms provide an alternative method to that described in the
previous subsection for solving the generalised word problem for $H$ in $G$.
The given word in $w \in A^*$ is reduced (in quadratic time, using the coset
multiplier automata) to the unique word $w^\prime$ in the language of the
coset word-acceptor for which $H{\overline w} = H{\overline {w^\prime}}$. 
Then $w \in H$ if and only if $w^\prime$ is the empty word.
The two methods of solving the generalised word problem are to some extent
complementary to each other, since there can exist $L$--quasiconvex subgroups
that are not coset automatic and vice versa, although the two concepts are
equivalent in word-hyperbolic groups.

The additional hypothesis required for the algorithm developed by Hurt
is the following generalisation of the fellow-traveller condition.
Let $\{W,M_a\}$ be the shortlex automatic coset system for
$(G,H)$ that we are trying to compute. Then, if
$(w^+,x^+) \in L(M_a)$ for some $a \in A$, there exists $h \in H$
such that $wa =_G hx$. The hypothesis is that there exists a constant
$k \geq 0$ such that for all such $w, x, a$ and $h$, and all $t \geq 0$,
we have $l_A( 
{\overline {w(t)}}^{-1} h {\overline {x(t)}} ) \leq k$.
In particular, taking $t=0$, we get $l_A(h) \leq k$, and so
in all such equations, only a finite number of elements $h$ occur.

One step in the algorithm is to define the states of the $M_a$ as
triples $(s_1,s_2,g)$, as in the automatic group case, but now
the initial states are $(s_0,s_0,h)$, where $s_0$ is the initial state
of $W$, and $h$ is one of the elements of $H$ occurring in the above
equations. So the $M_a$ are in fact constructed initially as
non-deterministic automata with multiple initial states,

If the hypothesis holds, then we shall say that $G$ is {\em strongly
coset automatic} with respect to $H$, and call $\{W,M_a\}$ a
{\em strong automatic coset system} for $(G,H)$.
It is proved in~\cite{Hurt} that word-hyperbolic groups are always
strongly coset automatic with respect to their quasiconvex subgroups.
It is easy to construct examples in which the hypothesis does not hold, by
choosing $H$ to be normal in $G$, in which case $G$ coset automatic with
respect to $H$ is equivalent to $G/H$ automatic, but we do not know of any
example in which ${\rm Core}_G(H) = 1$.

\section{HNN extensions and an example}
For the application to be described in this section, we need
to strengthen the hypothesis defined at the end of the preceding section
for strong automatic coset systems.

\begin{defn}
Let $\{W,M_a\}$ be a strong automatic coset system for $(G,H)$
with respect to the generating set $X$ of $G$. Let $Y$ be a finite
set of generators of $H$, and let $B = Y \cup Y^{-1}$.
Then $Y$ is said to be {\rm efficient} with respect to $\{W,M_a\}$ if, for any
$w,x \in L(W)$ and any $b \in  B, h \in H$ such that $wb =_G hx$, we
have either $h = 1$ or $h \in B$.
\end{defn}

We are not currently aware of any particular situations under which an
efficient generating set could be shown to exist; it would be
interesting to investigate this question. In specific examples of
automatic coset systems that we have calculated with the programs, it
is often possible to observe directly from the calculation that
a particular $Y$ is efficient.
The concept is useful to us here,
because it enables us to prove the following result about HNN extensions,
which can then be applied to a specific example. Note that a rather different
condition under which an HNN extension of an automatic group is 
asynchronously automatic has been proved by Shapiro in~\cite{Shap}, and
results of a similar nature
for amalgamated free products are proved in~\cite{BGSS}.

\begin{theorem}
Let $\{W,M_a\}$ be a strong automatic coset system for $(G,H)$,
let $G = \langle X\,|\,R \rangle$ be a finite presentation of $G$,
and suppose that $H$ has the efficient generating set $Y$.
Suppose also that $H$ is automatic, and 
let $\alpha$ be an automorphism of $H$ such that $\alpha(Y) = Y$.

Then the HNN extension
$$K = \langle X, z\,|\,R,\,z^{-1}yz = \alpha(y)\ (y \in Y) \rangle$$
is asynchronously automatic.
\end{theorem}
\begin{proof}
Let $T$ be a right transversal for $H$ in $G$. Then by the normal form
theorem for HNN extensions (see, for example, Theorem 2.1 (II), page 182
of~\cite{LS}), each element of $g \in K$ has a unique expression of the
form $$k = ht_1 z^{n_1} t_2 z^{n_2} \ldots t_r z^{n_r},$$
where $h \in H$, $t_i \in T$, $n_i \in \Z$, $t_i \not\in H$ for $i>1$
and $n_i \neq 0$ for $i < r$.

We use this normal form in the natural manner to construct a regular language
$L_K$ for $K$  on the alphabet $A \cup B \cup \{z^{\pm 1}\}$ where, as before,
$A = X \cup X^{-1}$ and $B = Y \cup Y^{-1}$. We are assuming that $H$
is automatic, so we can use the language $L_H$ of the word-acceptor from an
associated automatic structure with alphabet $B$ to obtain a word $w_h \in
L_H$ for the element $h \in H$ in the normal form.  For $T$ we choose the
image in $G$ of $L(W)$, and to represent $t_i$, we choose the unique
word $w_i \in L(W)$ with ${\overline {w_i}} = t_i$. This clearly yields a
regular language $L_K$ mapping bijectively onto $K$.

We now have to show how to construct the asynchronous multiplier
automata $M_c$ for $c \in A \cup B \cup \{z^{\pm 1}\}$.
Since this is fairly routine, we describe the construction in outline
only.
Suppose that $u, v \in L_K$ and $uc =_K v$, and let the
HNN normal form of $k = {\overline {u}}$
be $ht_1 z^{n_1} t_2 z^{n_2} \ldots t_r z^{n_r},$ as above.
If $c = z$ or $z^{-1}$, then the HNN normal form
for $kc$ in $K$ is just $ht_1 z^{n_1}  \ldots t_r z^{n_r \pm 1}$, and it is
easy to construct $M_c$.
So suppose $c \in A \cup B$.
We shall suppose that $n_r \neq 0$ and omit the details of the case
$n_r=0$, which are similar.
There exist words $c_1 \in B^*$ and $c_2 \in L(W)$ such that $c =_G c_1c_2$.
Let $l_B(c_1) = k$.
Then, from the assumptions that the generating
set $Y$ of $H$ is efficient and that $\alpha(Y) = Y$, it follows that the
HNN normal form in $K$ for $kc$ is
$$kc = h^\prime t_1^\prime z^{n_1} t_2^\prime z^{n_2} \ldots t_r^\prime
z^{n_r}{\overline {c_2}},$$
where there are elements $x_i, y_i \in H\ (1 \leq i \leq r)$,
all having $B$--length at most $k$,
such that $z^{n_r} {\overline{c_1}} = y_r z^{n_r}$,
$t_i y_i = x_i t_i^\prime$ for $1 \leq i \leq r$,
$z^{n_i} x_{i+1} = y_i z^{n_i}$ for $1 \leq i < r$,
and $hx_1 = h^\prime$.
Thus we have
$u = w_hw_1 z^{n_1} \ldots w_r z^{n_r}$ and
$v = w_{h^\prime}w_1^\prime z^{n_1} \ldots w_r^\prime z^{n_r} c_2$,
where $w_h, w_{h^\prime} \in L_H$  map onto $h, h^\prime \in H$,
and $w_i, w_i^\prime \in L(W)$ map onto $t_i, t_i^\prime \in T$ for
$1 \leq i \leq r$.

The multiplier $M_c$ proceeds by reading the words $w_h$ and
$w_{h^\prime}$ in parallel at the same rate, then the $z^{n_1}$
together, then $t_1$ and $t_1^\prime$ together, and so on. If either
of $w_h$ or $w_{h^\prime}$ is longer than the other, then it will
wait at the end of the shorter one until the longer word has been read,
and similarly for $t_i$ and $t_i^\prime$. (This explains
why $M_c$ needs to be asynchronous. Although $|l(w_h) - l(w_{h^\prime})|$
and $|l(t_i) - l(t_i^\prime)|$ are all bounded, there is no bound on $r$,
and so one of the two tapes of the input of $M_c$ may conceivably
get indefinitely ahead of the other; indeed, we have verified that this
really can happen in the example below.)

Of course, if either of the two words input to $M_c$ is not in $L_K$, or
if they do not both have the same pattern with respect to the
occurrences of $z$, then they are rejected. Otherwise, if after $t$
transitions, $M_c$ has read $\phi(t)$ symbols from $u$ and $\psi(t)$ from
$v$, then the element
$g(t) = {\overline {u(\phi(t))}}^{-1}{\overline {v(\psi(t))}}$ of $K$
is remembered as a function of the state of $M_c$. As in the synchronous case,
it is sufficient to show that $l(g(t))$ is bounded.

There are four essentially different situations that occur as the words
$u,v$ are read.
\begin{enumerate}
\item[(i)] $u(\phi(t))$ and $v(\psi(t))$ are prefixes of $w_h$ and
$_{h^\prime}$, where $|\phi(t) - \psi(t)|$ is bounded.
Then the the boundedness of $l(g(t))$ from the automaticity of $H$,
and the fact that $hx_1=h^\prime$ with $l(x_1) \leq k$.
\item[(ii)] $u(\phi(t)) = w_hw_1 z^{n_1} \ldots w_i(s_1)$ for some $i$ and
some prefix $w_i(s_1)$ of $w_i$, and
$v(\psi(t)) = w_{h^\prime}w_1^\prime z^{n_1} \ldots w_i^\prime(s_2)$,
where $|s_1-s_2|$ is bounded.
Then $g(t) =
{\overline {w_i(s_1)}}^{-1} x_i {\overline {w_i^\prime(s_2)}}$,
and its boundedness follows from the assumptions that $l_B(x_i) \leq k$
and that $\{W,M_a\}$ is a strong automatic coset system for $(G,H)$.
\item[(iii)] $u(\phi(t)) = w_hw_1 z^{n_1} \ldots w_iz^{m_1}$ for some $i$
and some $m_1 \leq n_i$, and
$v(\psi(t)) = w_{h^\prime}w_1^\prime z^{n_1} \ldots w_i^\prime z^{m_2}$,
where $|m_1-m_2| \leq 1$.
Then $g(t)
z^{-m_1}y_iz^{m_2}$, and its boundedness follows from $l_B(y_i) \leq k$
and the assumption that $\alpha(Y) = Y$.
\item[(iv)]  $\phi(t) > l(u)$ and $\psi(t) \geq l(v) - l(c_2)$.
Then $l(g(t)) \leq l(c_2)$ which is clearly bounded.
\end{enumerate}

This completes the proof of the theorem.
\end{proof}
\vspace{1mm}

As an application, we shall use this theorem together with the results of some
machine computations that were done with \kbmag, to prove that the group
defined by the presentation
\begin{eqnarray*}
\lefteqn{\langle \,a,b,r,t,x,z \, | \,}\\
&&xaxa=t, bxbx=t, bbtaa=t, a^{-1}br=ra^{-1}b, zt=tz, btaz=zbta \, \rangle
\end{eqnarray*}
is asynchronously automatic.

This group, which we shall denote by $K$, was originally proposed by Mark Sapir
as a possible building block in his attempts to construct groups with
given Dehn functions. However, he later found a different approach to
his problem, and so the example is no longer relevant from that viewpoint.
He had hoped that it could be proven automatic, but the methods we have been
discussing in this paper only appear to be sufficient to prove it
asynchronously automatic.

The computer programs could make no progress with the presentation as given
above, but matters improved after manipulating it a little.
Eliminating $t=bxbx$, we get
\begin{eqnarray*}
\lefteqn{\langle \,a,b,r,x,z \, | \,xaxa=bxbx, bbxbxaa=xbx, }\\
&&a^{-1}br=ra^{-1}b, zbxbx=bxbxz, bbxbxaz=zbbxbxa
\, \rangle.
\end{eqnarray*}
Now, putting $u=xa$ and $v=bx$, and eliminating $a=x^{-1}u=v^{-1}bu$ and
$x = b^{-1}v$, we get
\begin{eqnarray*}
\lefteqn{\langle \,u,v,b,r,z \, | \,u^2=v^2, bvbuv^{-1}bu=b^{-1}v^2, }\\
&&u^{-1}b^{-1}vbr=ru^{-1}b^{-1}vb,
zv^2=v^2z, bvbuz=zbvbu
\, \rangle.
\end{eqnarray*}
Finally, using $u^2=v^2$ to simplify the second relation , we get
\begin{eqnarray*}
\lefteqn{\langle \,u,v,b,r,z \, | \,u^2=v^2, bvbu=b^{-1}ub^{-1}v,}\\
&& u^{-1}b^{-1}vbr=ru^{-1}b^{-1}vb, zu^2=u^2z, zbvbu=bvbuz
\, \rangle,
\end{eqnarray*}

This is now visibly an HNN extension of the group
$$G = \langle \,u,v,b,r \, | \,u^2=v^2, bvbu=b^{-1}ub^{-1}v,
u^{-1}b^{-1}vbr=ru^{-1}b^{-1}vb \, \rangle.$$
with respect to the subgroup $H = \langle\,u^2, bvbu\,\rangle$, where
$H$ is centralised by the new generator $z$. (In fact $G$ is itself
an HNN extension with extra generator $r$, but we shall not make use
of that fact.)

Running the automatic coset system program from \kbmag\ on the subgroup
$H$ of $G$ verifies that $G$ is strongly coset automatic with respect to $H$.
(The coset word acceptor has 302 states, and the coset multiplers
about 1400 states.) The presentation of $H$ computed by the program
proves that $H$ is free of rank 2, and so it is certainly automatic.
The programs can also be used to verify that
the set $Y = \{u^2, bvbu^{-1} \}$ is an efficient generating set
for $H$. (Briefly, this is done by constructing the multiple
initial state multiplier automata for the elements
$u^2$ and $bvbu^{-1}$. The elements of $H$ corresponding to the
initial states of these automata can then be inspected from the output,
and it turns out that these are just the identity and elements of
$B = Y \cup Y^{-1}$.) We can now deduce from the theorem that
Sapir's group $K$ is asynchronously automatic.

As a final remark about this example, it turns out (again using
calculations carried out by \kbmag) that the subgroup $H$
is not $L$--quasiconvex, where $L$ is the language of the word-acceptor
of the shortlex automatic structure of $G$. The element
$(bub^{-1}v^{-1})^n(b^{-1}vbu^{-1})^n$ of $L$ lies in $H$ for all
$n \geq 0$, but the coset representative of $(bub^{-1}v^{-1})^n$
in the language of the coset word acceptor is $b^{2n}$.

\Addresses\recd

\end{document}